\documentstyle[12pt]{article}
\makeatletter
\newcommand{\singlespacing}{\let\CS=\@currsize\renewcommand{\baselinestretch}{1}\normalsize\CS}
\newcommand{\oneandahalfspacing}{\let\CS=\@currsize\renewcommand{\baselinestretch}{1.125}\tiny\CS}
\newcommand{\doublespacing}{\let\CS=\@currsize\renewcommand{\baselinestretch}{1.5}\tiny\CS}
\begin{document}
\thispagestyle{empty}
\singlespacing
\parskip = 10pt
\newcommand{\namelistlabel}[1]{\mbox{#1}\hfil}
\newenvironment{namelist}[1]{%
\begin{list}{}
{
\let\makelabel\namelistlabel
\settowidth{\labelwidth}{#1}
\setlength{\leftmargin}{1.1\labelwidth}
}
}{%
\end{list}}
\newtheorem{theorem}{\bf Theorem}[thetheorem]
\newtheorem{corollary}{\bf Corollary}[thecorollary]
\newtheorem{remark}{\bf Remark}[theremark]
\newtheorem{lemma}{\bf Lemma}[thelemma]
\def\theequation{\thesection.\arabic{equation}}
\def\thetheorem{\thesection.\arabic{theorem}}
\def\theremark{\thesection.\arabic{remark}}
\def\thecorollary{\thesection.\arabic{corollary}}
\def\thelemma{\thesection.\arabic{lemma}}
\newcommand{\bsp}{\begin{sloppypar}}
\newcommand{\esp}{\end{sloppypar}}
\newcommand{\be}{\begin{equation}}
\newcommand{\ee}{\end{equation}}
\newcommand{\beanno}{\begin{eqnarray*}}
\newcommand{\inp}[2]{\left( {#1} ,\,{#2} \right)}
\newcommand{\dip}[2]{\left< {#1} ,\,{#2} \right>}
\newcommand{\disn}[1]{\|{#1}\|_h}
\newcommand{\pax}[1]{\frac{\partial{#1}}{\partial x}}
\newcommand{\tpar}[1]{\frac{\partial{#1}}{\partial t}}
\newcommand{\xpax}[2]{\frac{\partial^{#1}{#2}}{\partial x^{#1}}}
\newcommand{\pat}[2]{\frac{\partial^{#1}{#2}}{\partial t^{#1}}}
\newcommand{\ntpa}[2]{{\|\frac{\partial{#1}}{\partial t}\|}_{#2}}
\newcommand{\xpat}[2]{\frac{\partial^{#1}{#2}}{\partial t \partial x}}
\newcommand{\npat}[3]{{\|\frac{\partial^{#1}{#2}}{\partial t^{#1}}\|}_{#3}}
\newcommand{\lnpat}[4]{{\|\frac{\partial^{#1}{#2}}{\partial
t^{#1}}\|}_{L^{#3}(H^{#4})}}
\newcommand{\xkpat}[3]{\frac{\partial^{#1}{#2}}{\partial t^{#3} \partial x}}
\newcommand{\jxpat}[3]{\frac{\partial^{#1}{#2}}{\partial t \partial x^{#3} }}
\newcommand{\eeanno}{\end{eqnarray*}}
\newcommand{\bea}{\begin{eqnarray}}
\newcommand{\eea}{\end{eqnarray}}
\newcommand{\ba}{\begin{array}}
\newcommand{\ea}{\end{array}}
\newcommand{\nno}{\nonumber}
\newcommand{\dou}{\partial}
\newcommand{\bc}{\begin{center}}
\newcommand{\ec}{\end{center}}
\newcommand{\bb}{\mbox{\hspace{.25cm}}}
\newcommand{\re}{{\bf |\!\!R}}
\newcommand{\call}{{\cal L}}
\newcommand{\calla}{{\cal L}^*}
\newcommand{\caq}{{\cal Q}_h}
\newcommand{\calm}{{\cal M}(u)}
\newcommand{\calma}{{\cal M^*}(u)}
\newcommand{\dm}{\displaystyle}
\newcommand{\tilf}{\tilde{f_{uu}}}
\title{\bf  Substantial Riemannian Submersions of $S^{15}$
with 7-dimensional Fibres}
\author{Akhil Ranjan}
\maketitle
\thispagestyle{empty}
\begin{abstract}
In this paper we show that a substantial Riemannian submersion of $S^{15}$
with 7- dimensional fibres is congruent to the standard Hopf fibration.
As a consequence we prove a slightly weak form of the Diameter Rigidity
theorem for the Cayley plane which is considerably stronger than the very
recent Radius Rigidity Theorem of Wilhelm.
\end{abstract}
\setcounter{equation}{0}
\thispagestyle{empty}
\section{\bf Introduction}
\label{sec-intro}
\setcounter{subsection}{0}
\setcounter{equation}{0}
The study of Riemannian foliations of Euclidean spheres
with no restrictions on the the geometry of leaves was initiated in
\cite{GG3} and a complete classification was obtained for all the
foliations having leaves of dimensions less than or equal to
three.  Earlier, the classification under the very strong
hypothesis of the leaves being totally geodesic was accomplished
in \cite{E},\cite{R} and \cite{GWZ}. A remarkable discovery in \cite{GG3} was
the
fact that for substantial or highly nonflat (see \cite{GG3} for
definition) Riemannian foliations of $S^{n}$, the horizontal
holonomy group reduces to a compact Lie group.  The low
dimensional Riemannian foliations could then be shown to be the
orbits of this group.   Thus proving substantiality was a
crucial step and it was carried out for the leaves of dimensions
less than or equal to three.  As a result it is known in
particular that if $\pi: S^n \rightarrow M$ is a Riemannian
submersion with connected fibres, then either it is congruent to
a Hopf fibration or we have $\pi:S^{15}\rightarrow M^8$.
Now the reduction to a compact Lie group is valid for
substatial foliations without restrictions on leaf dimension.
Hence it is natural to ask the following question:

{\it If the Riemannian submersion $\pi :S^{15}\rightarrow M^8$
is substantial along every leaf, is it congruent to the Hopf
fibration?}\\
Since the question is algebraic in nature ,this author was
under the impression that it must have been settled by now.
However looking at a recent preprint \cite{W} it was realised that
such wasnot the case. In \cite{W} the following lemma occurs

 {\it Let $\pi:S^{15}(1)\rightarrow V$ be a Riemannian submersion
with connected 7-dimensional fibres,and let G be the set of
points $v \in V$ so that $\pi^{-1}(v)$ is totally geodesic. Then
either G is discrete or G is totally geodesic and an
isometrically embedded copy of $S^l(1/2)$ for some $1\leq l\leq 8$. }

  Now in \cite{GG3} itself it is implicit that if $G\neq\phi$,then
$\pi$ will be substantial along every leaf and it should be the
case that G is all of V.  This would simplify much of the
subsequent work in \cite{W}.  Therefore, it definitely seems worthwhile
to set the record straight in this matter.
 In this paper we prove the following:
\begin{theorem}
If ${\bf \pi:S^{15}\rightarrow M^8}$ is a Riemannian
submersion with connected fibres which is substantial along each
leaf, then it is congruent to the Hopf fibration.
\end{theorem}
  As a corollary we also prove
\begin{theorem}
{\bf(Weak Diameter Rigidity for $CaP^2$)} Let $ M$
have sectional curvatures $\geq 1,\,diam(M)=\pi/2$.  Moreover,let M
admit an equilateral triangle of sides $\pi/2$ and have the
integral cohomology ring same as that of Cayley plane. Then M is
isometric to the standard $CaP^2$ with $ 1\leq K\leq 4$.
\end{theorem}
This almost ties up the loose end in \cite{GG2} and both the
$Radius\, Rigidity\, Theorem$ in \cite{W} and the
$corollary II$ in \cite{D} are considerably improved.\\
{\bf Acknowledgement:} The methods in this paper are direct
generalisations of those developed in \cite{GG3} and the author
acknowledges this without hesitation.  Also the author would like to
thank K. Ramachandran for bringing Wilhelm's paper \cite{W} to his notice
and Prof.  A.K. Pani who spent valuable time to help me typeset
the manuscript on {\it LATEX}.

  \section{The Group of the Submersion} In this section we
determine the Lie group which could possibly occur as the structure
group (holonomy group) of the Riemannian submersion $\pi:S^{15} \rightarrow
M^8$
when we assume it to be substantial.
\begin{theorem}
   If $\pi: S^{15} \rightarrow M^8$ is a Riemannian submersion which is
substantial along each leaf, then the holonomy group G of $\pi$ is
such that LieG = so(8).  Hence G is either SO(8) or Spin(8).
\end{theorem}
{\bf Proof} : By \cite{GG3} ,the holonomy group G of $\pi$ is a compact
Lie group of dimension $\leq 28$ and it acts transitively on
the fibres. The fibres of $\pi$ are well known to be homotopy
7-spheres (see \cite{Br}). We conclude therefore that G can  only be one
of the following:
$$ SO(8),\,Spin(8),\,Spin(7),\,SU(4),\,U(4),\,{\rm and}\,Sp(2)$$.
  (see \cite{Be}, p. 195).

  Now by the long exact sequence of homotopy of $\pi$ (or directly
from \cite{Br}) $M^8$ is a homotopy 8-sphere and by the Thom-Gysin sequence
of $\pi$, the Euler class of this 7-sphere fibration is the generator of
$H^8(M,Z\!\!\!Z)\approx Z\!\!\!Z$.  This enables us to rule out  Spin(7),
since the Euler class of a Spin(7) bundle must be torsion (this is
because $H^*(BSpin(7),I\!\!\!\!Q) = I\!\!\!\!Q[p_1,p_2,p_3]$ is generated
by Pontrjagin classes only).  Since $Sp(2) \subset SU(4) \subset U(4)$,
if we rule out U(4) then the other two will be ruled out automatically.
Now $H^*(BU(4),Z\!\!\!Z)= Z\!\!\!Z[c_1,c_2,c_3,c_4]$ is generated by
Chern classes and $c_4$ also represents the Euler class $\chi$.  Since
M is a homotopy sphere of dimension 8, by Bott's Integrality Theorem
(see \cite{H}, p. 279) $c_4$ is divisible by 3! in $H^8(M,Z\!\!\!Z)$.  But then
$\chi$ can not be a generator - a contradiction!
Therefore, we conclude that $G= Spin(8)$ or $SO(8)$.        \,\,\hfill{\bf
q.e.d.}

\noindent
{\bf REMARK:} Since  there are no Whitney classes $w_1$ and $w_2$ we can fix
the group to be $Spin(8)$ which we do from  now on.
\begin {corollary}
  Let ${\cal X}(N)$ denote the Lie algebra of smooth vector fields on a
smooth manifold N. For $b\in  M$,the map
$$ A :\bigwedge\,\!\!^{2}(T_bM) \rightarrow {\cal X} (\pi^{-1}(b))$$
given by $ A(x\wedge y)= A_xy$ is injective and its image is a Lie subalgebra
isomorphic to so(8).
 \end{corollary}
{\bf proof}: See \cite{GG3}.
\section {The associated principal bundle}
  Having found the holonomy group $ \pi:S^{15}\rightarrow M^8$ we cannot
apply the methods of \cite{GG3} directly at this stage. In fact the standard
Hopf fibration $S^{15}\rightarrow S^8(1/2)$ does not occur via orbits
of any Lie group action on $S^{15}$. We therefore adopt the strategy
of looking at the associated principal $Spin(8)$ bundle
$$\tilde{\pi}\,:\,E\,\rightarrow\, M$$
and study the induced connection on it.
  To construct the associated principal bundle we first note that the
holonomy Lie algebra $so(8)\subset {\cal X} (\pi^{-1}(b))$ for each $b\in M$
gives an additional geometric structure on the fibres which is preserved under
the horizontal holonomy displacement (\cite{GG3}, lemma 2.12, and prop.2.13).
Moreover,the fibres being homogeneous under $Spin(8)$ action (though
not via isometries {\it apriori}) become diffeomorphic to $S^7(1)$.
 Now $Spin(8)$ can act transitively on $S^7$ in three ways: via
$\rho,\, \sigma_+,\,and\, \sigma_-$, where $\rho$ is the natural double cover
representsation in $I\!\!R^8$ and the other two are the two spin
representations of $Spin(8)$.  However, from the geometrical point of
view, the three are equivalent as all three are two sheeted covers
of $SO(8)$ and one can pass from one to any other using covering
transformations. These are naturally {\it outer} automorphisms
of $Spin(8)$. To construct the principal bundle we take any one of these
represntations, say $\rho $.  The resulting $Spin(8)$ action on $S^7$
gives rise to a copy of $so(8)$ inside ${\cal X}(S^7(1))$. Let for each
$b\in M,F_b$ denote the fibre $\pi^{-1}(b)$ and define $E_b$ to
be the set of all diffeomorphisms $\phi :S^7(1)\rightarrow F_b$ which
send this copy of $so(8)$ onto the holonomy Lie algebra sitting in
${\cal X}(F_b)$ isomorphically. Clearly $E_b$ is a copy of the group
Aut(so(8)) whose identity component is just $PSO(8)$.  We now let $E$ to
be the disjoint union of these $E_b$ as $b$ runs through the set $M$.
The set E has a natural smooth structure and an obvious smooth projection
$\tilde\pi$ onto M with fibre over any point $b$ being exactly $E_b$.
The group $Aut(so(8))$ acts on $E$ from the right freely and the orbits
are the fibres $E_b$.  We now do two things for the sake of definiteness:
(i) reduce the group to the identity component $PSO(8)$,this is possible
due to simple connectivity of M and (ii){\it reduce} the group further
to $Spin(8)$.  This is possible due to the absence of the first two
Whitney classes.  We continue to denote the reduced space by E only.
Thus E is now the associated principal Spin(8) bundle we were looking for.
If Spin(7) is the isotropy group at some chosen base point of our model
sphere $S^7(1)$, then $E/Spin(7)$ is our $S^{15}$ that we started with and
we have a tower of smooth fibrations
$$E\, \rightarrow\, S^{15}\, \rightarrow\, M$$
  The composite is the map $\tilde{\pi}$.
\section{The Riemannian structure on E}
In this section we describe a
metric on E which comes naturally due to the geometric considerations
and which makes the above tower,a tower of Riemannian submersions. This is
done in several stages.

First we note that for each $\phi\in E$,there is a
decomposition of the tangent space into two complemntary subspaces.  The
first of  these is the tangent space to the orbit of Spin(8) through $\phi$
and the other is the space generated by the holonomy of $\pi$.  To elaborate
this, let $\phi\in E_b$ and $\gamma $ a path starting from $b$.  Then we get
a one parameter family of diffeomorphisms $\tau_{\gamma(t)}$ from $F_b$ to
$F_{\gamma(t)}$ generated by the holonomy displacement along $\gamma$.
This gives a path $\tau_{\gamma(t)}o\phi$ in E. Its derivative at $t=0$ gives
a member of the second complementary subspace.  As $\gamma$ varies the
full subspace is obtained. We naturally declare the break-up as orthogonal
and metrize the second part by the inner product on $T_bM$ to which it
projects isomorphically under $\tilde{\pi}_*$.  Henceforth we will denote
the second subspace $\cal{H}_\phi$ and call it the horizotal space
(relative to $\tilde \pi$). We remark that this collection of horizontal
spaces gives exactly the connection on the principal bundle which is
associated to the connection given by the horizontal spaces in $S^{15}$.
We also denote the space complementary to $\cal{H}_\phi$ by $\cal{V}_\phi$.
It remains to define a metric on this {\it vertical} part. This is more
delicate.  We proceed as follows:
 We already have a copy of $so(7)$ inside $so(8)$ as the isotropy Lie
algebra of a chosen base point in our model $S^7$. This naturally breaks
$so(8)$ in a unique manner as $$so(7)\oplus I\!\!R^7$$ each part being
an irreducible $so(7)$ module. This is also the Cartan decomposition
(see \cite{He}).
  Corresponding to each $X\in so(8)$,
we have a vertical vector field $\bar{X}$ on $E$ (see \cite{BC},p.39).
Thus we also have a decomposition
$$\cal{V}_\phi =\cal{V}_\phi'\oplus\cal{V}_\phi''$$
  This too we declare to be orthogonal. On the first part we put the metric
coming from the bi-invariant innerproduct of $so(7)$ and on the second
part the metric pulled back from $S^{15}$ to which it goes injectively
under $\bar{\pi_*}$.  Note that its image is precisely the vertical space
 at $\bar{\pi}(\phi)$ in $S^{15}$. This completes the description of
the Riemannian structure on $E$. With this metric on $E$,we have a tower of
Riemannian submersions
$$E\rightarrow S^{15}\rightarrow M^8$$
 The first one is a $Spin(7)$ principal bundle with totally geodesic
fibres,each isometric to  $Spin(7)$ with the Cartan-Killing metric
while the composite is the principal $Spin(8)$ bundle whose fibres
{\it apriori} have varying metrics. The decomposition
$$TE=\cal{V}\oplus\cal{H}$$
is preserved under the $Spin(8)$ action while the finer decomposition
$$TE=\cal{V}'\oplus\cal{V}''\oplus\cal{H}$$
is preserved under the $Spin(7)$ action.
Moreover, the $Spin(7)$ action is evidently via isometries.
\section{The Metric on the Fibres of $\tilde{\pi}$}
 In this section we will see that though the metric on the fibres
$E_b, b \in M$ could be different,the variation is very limited in
nature. More precisely we have the following:
\begin{theorem}
The fibres $E_b, b \in M$ of $\tilde\pi$ are isometric
to $Spin(8)$ furnished with a left invariant metric which is also
right invariant under $Spin(7)$ action.
\end{theorem}
To prove this result we need to analyse the integrability tensor
of $\tilde\pi$.  So let $\tilde{A}$ denote the integrability tensor of
$\tilde\pi$ and $\bar{A}$ that of $\bar{\pi}$.
\begin{lemma}
Let $x,y\in T_bM$. The vertial field $\tilde{A}_xy$
along the fibre $E_b$ is $\bar{Z}$,for some $Z \in so(8)$. Further,if
$Z = Z_1 + Z_2$ corresponding to $so(8) = so(7) + I\!\!R^7$,then
$\bar{Z_1} = \bar{A}_xy\, and\, \bar{Z_2} = A_xy$
\end{lemma}
{\bf Proof} $\tilde{A}_xy = 1/2[X,Y]^v$ is just the curvature of the
cnnection on the $Spin(8)$ principal bundle alluded to earlier and
hence clearly comes from a suitable element of its Lie algebra.
Also  $[X,Y]^v = [X,Y]^{v'}+[X,Y]^{v''} =2(\bar{A}_xy+A_xy)$ is also
$2(\bar{Z}_1 + \bar{Z}_2)$.  By uniqueness of decomposition into components,
the result follows.  \hfill{\bf q.e.d.}
\begin{corollary}
For $x,y\in T_bM,\tilde{A}_xy$ is of constant length
along $E_b$.
\end{corollary}
{\bf Proof}:  $\bar{Z}_1$ is of constant length since it comes from an
$so(7)$ element and $\bar{Z}_2 = A_xy$ is of constant norm since it
is the basic lift of the correspondig field along $F_b\subset S^{15}$.
That it is of constant norm there is well known (see \cite{GG3}). Since the
two components are also mututually orthogonal everywhere the corollary
follows.                                                     \hfill{\bf q.e.d.}
\begin{lemma}
There is an {\it orthonormal framing} of $E_b$ which is
generated by a basis of $so(8)$.
\end{lemma}
{\bf Proof}: Choose a basis $\{x_i:0\leq i\leq 7\}$ of $T_b M$ so that
$\{A_{x_0}x_i:1\leq i\leq 7\}$ forms an orthonormal framing of $F_b$.
 Set  $v_i = A_{x_0}x_i$, and consider the fields
 $$\{v_i:1\leq i\leq 7\}\cup \{\bar{A}_{v_j}v_k : 1\leq j<k\leq 7\}$$
That the first set of 7 vetor fields is an orthonormal frame along $E_b$ is
obvious. That it is generaated by $so(8)$ elements is one of the contention
of the lemma 5.1 above. As for the second set of 21 fields
we note that if $v_i =\bar{X}_i$ for suitable $X_i$ in the $I\!\!R^7$
component of $so(8)$ then $\bar{A}_{v_j}v_k$ is generated by
$[X_j,X_k]$.  There is no need to project to the $so(7)$ part since it
is already there.  (This is a well known property of Cartan decomposition).
 That these have pairwise constant inner-products now follows from the
very way the metric was defined on the $\cal{V}''$ part of $TE$.
The two sets are cleary mutually orthogonal.  We also note that
$\,\{[X_i,X_j],1\leq i<j\leq 7\}\, $ is a basis of $so(7)$.  An application
of Gram-Schmidt orthonormalisation on the second set now gives the
required framing.                  \hfill{\bf q.e.d.}\\
{\bf proof of the theorem}: $E_b$ can  be identified with $Spin(8)$ after
choosing some base point on it. The left invariant vector fields of $Spin(8)$
are mapped to the $I\!\!R$-span of the above mentioned basis.  Left invariance
of the metric of $E_b$ now follows exactly as in \cite{GG3}.  Right invariance
under $Spin(7)$ is from construction.                          \hfill{\bf
q.e.d.}  \\
{\bf Caution:}  Left and right in our situation are opposite to those
in \cite{GG3}.
\begin{corollary}
Each fibre $F_b$ of $\pi$ is a sphere of constant
sectional curvature.
\end{corollary}
{\bf proof}: Each $F_b$ is isometric to the quotient of $E_b$ under the
group $Spin(7)$ of isometries of $E_b$ acting from the right. Hence
$Spin(8)$ acts on $F_b$ via isometries from the left making it a
homogeneous Riemannian space. Since the isotropy group $Spin(7)$ acts
transitively on tangent two-planes the sectional curvatures are constant
pointwise.  It follows they are the same constant everywhere (see \cite{KN}).
\hfill{\bf q.e.d.}
\section{Left invariant Metrics on $Spin(8)$ invariant under
$AdSpin(7)$}
Let $$so(8) = so(7)\oplus I\!\!R^7$$ be the Cartan decomposition which
we note is also an $so(7)$-module decomposition into its {\it isotypical}
components.  With this notation we have the following:
\begin{theorem}
Let $g$ denote the Cartan Killing metric on $so(8)$.
Then $g = g_1\oplus g_2$ corresponding to the Cartan decomposition
and any $Ad Spin(7)$ invariant metric is of the type
$$g = g_1\oplus cg_2$$ (upto an overall scalar multiple).
\end{theorem}
{\bf proof} The first claim is well known (see \cite{He}). The second claim
follows since under an $Ad Spin(7)$-invariant metric, $so(7)^{\perp}$
must be an $so(7)$ module.  This forces the above decomposition to be
orthogonal.  Further any $Spin(7)$ invariant metric on the natural
module $I\!\!R^7$ is a scalar multiple of the standard Euclidean
metric.                                                \hfill{\bf q.e.d.}
\begin{corollary}
There is a smooth positive real valued function
$c$ on $M$ such that $E_b$ is isometric to $Spin(8)$ with the left
invariant metric $g_b = g_1\oplus c(b)g_2$.
\end{corollary}
{\bf proof}: Just observe that there is no change of scale in the
$\cal{V}'$ part of the tangent space of $E_b$ as the orbits of the
subgroup $Spin(7)$ are all isometric to each other.

\hfill{\bf q.e.d.}
\begin{corollary} For any smooth path in $M$ the holonomy displacement
of fibres in $E$ preserves the $\cal{V} =\cal{V}'\oplus \cal{V}''$
decomposition.  Moreover,it is an isometry on the first part and a
dilatation on the second part.
\end{corollary}
{\bf proof} Let $\gamma$ be a smooth path in $M$ starting from $b_1$
and ending at $b_2$.  Let $\phi _1\in E_{b_1}$ and $\phi _2 =\tau (\phi _1)$
be in $E_{b_2}$.  On identifying these fibres with $Spin(8)$, it is
clear that the holonomy diffeomorphism gets identified with $L_g$
for a suitable $g\in Spin(8)$. Since $g_{b_{2}}$, the metric on $E_{b_2}$
pulls back under $L_g$ to a left invariant metric which is also right
invariant under $Spin(7)$, so it must be of the form $g_1\,\oplus\,
\frac{c(b_2)}{c(b_1)}\,.\,g_2 $ ( see corollary 6.4 for notation ).\hfill{\bf
q.e.d.}

\section{Proof of the Main Theorem}

In this section we complete the proof of Theorm1. 1 stated in the
introduction. We observe that a consequence of the last corollary is
that the holonomy displacement in $S^{15}$ is via isometries upto
scalings factor. This forces  all the operators $T^x$ on vertical
vectors arising out of the {\it second fundamental form} of the
fibres to be scalar multiples of identity. This in turn forces the
largest and the smallest fibres to be great spheres due to
isoparametricity. Hence all the fibres are great spheres.  Appealing
now to the classification of parallel great sphere fibrations
as in \cite{R} or \cite{GWZ} we conclude that $\pi$ is a Hopf fibration.
\hfill{\bf q.e.d.}

\begin{corollary}
{\bf (Generalized Wilhelm's Lemma)}:  Let
$$\pi :S^{15}(1)\rightarrow M^8$$ be a Riemannian submersion with
connected 7-dimensional fibres and let
$$G\,=\,\{b\in M:\pi^{-1}(b)\,\, {\rm \,is\,\, a\,\, great\,\, sphere\,\,} \}$$
then either $G$ is empty or all of $M$.
\end{corollary}
{\bf proof}: If $G$ is nonempty let $b\in G$.  For any $p\in F_b$
and any $x\in {\cal H}_p$, the linear map
$$A_x :{\cal H}_p \rightarrow {\cal V}_p$$
is surjective (see lemma 4.2, \cite{GG3}).  It follows that $A_{c_x'(t)}$
is surjective for all $t$, where $c_x$ denotes the geodesic with
initial vector $x\in {\cal H}_p$ (see the statement preceding the
corollary 2.9 in \cite{GG3}).  Hence $\pi$ is substantial along every
fibre.                                               \hfill{\bf q.e.d.}
\section{Weak Diameter Rigidity for $CaP^2$}
{\bf proof}: Let $x$,$y$,and $z$ be mutually distance $\pi /2$
apart.  In this situation $y$ and $z$ both are in the dual set $\{x\}'$
of $x$.  (See \cite{GG2} and \cite{W} for more details about dual sets. )
{}From \cite{GG2} we know that there is a Riemannian submersion
$$exp_x :S_x\,\rightarrow\,\{x\}'$$
and similarly for $y$ and $z$.  As argued in \cite{W} this forces at
least one fibre to be totally geodesic in each case.  But then $exp_p$
is congruent to the Hopf fibration and this makes the space isometric
to the standard $CaP^2$ as proved in \cite{GG2} and \cite{W}.
\hfill{\bf q.e.d.}\\
{\bf Remark}: This rigidity is clearly stronger than even
{\it corollary II} of \cite{D} where each pair of points seperated by
a distance of $\pi /2$ is required to be completed into an equilateral
triangle. \\
 It is clearly an interesting problem to find reasonable conditions
which will force substantiality.  One such as we saw is total geodeity
of a single fibre.  The author proposes to discuss this elsewhere.

\noindent
Deparment of Mathematics\\
Indian Institute of Technology, Bombay\\
Powai, Bombay - 400 076, INDIA.\\
email: aranjan@ganit.math.iitb.ernet.in
\end{document}